\documentclass[12pt]{amsart}
\usepackage{amsmath,amssymb,amsfonts,amscd,verbatim,graphicx}
\usepackage{xypic}
\newtheorem{theorem}{Theorem}[section]
\newtheorem{lemma}[theorem]{Lemma}

\theoremstyle{definition}
\newtheorem{definition}[theorem]{Definition}

\theoremstyle{remark}

\numberwithin{equation}{section}

\theoremstyle{definition}

\newcommand{\co}{\colon\thinspace}

\newcommand{\nl}{\hfil\break}

\begin{document}

\title{Graphs, links, and duality on surfaces}
\author{Vyacheslav Krushkal}
\address{Department of Mathematics, University of Virginia, Charlottesville, VA 22904}
\email{krushkal\char 64 virginia.edu}

\thanks{Supported in part by NSF grants DMS-0605280 and DMS-0729032}

\begin{abstract}
We introduce a polynomial invariant of graphs on surfaces, $P_G$, generalizing the classical Tutte polynomial.
Topological duality on surfaces gives rise to a natural duality result for $P_G$, analogous to the
duality for the Tutte polynomial of planar graphs. This property is important from the perspective of
statistical mechanics, where the Tutte polynomial is known as the partition function of the Potts model.
For ribbon graphs, $P_G$ specializes to the well-known Bollob\'{a}s-Riordan
polynomial, and in fact the two polynomials carry equivalent
information in this context. Duality is also established for a multivariate version of the polynomial $P_G$.
We then consider a $2$-variable version of the Jones polynomial for links in thickened surfaces, taking
into account homological information on the surface.
An analogue of Thistlethwaite's theorem is established for these generalized Jones and Tutte polynomials
for virtual links.
\end{abstract}

\maketitle

\section{Introduction}

The Tutte polynomial $T_G(X,Y)$ is a classical invariant in graph theory (see \cite{T1, T2, B}), reflecting
many important combinatorial properties of a graph $G$. For example, the chromatic polynomial, whose values at positive
integer values of the parameter $Q$ correspond to the number of colorings of $G$ with $Q$ colors, is a one-variable
specialiazation of $T_G$. The Tutte polynomial is also important in statistical mechanics, where it arises as the
partition function of the Potts model, cf \cite{S}.

Two properties of the Tutte polynomial are particularly important in these contexts: the contraction-deletion rule, and the duality
\begin{equation} \label{duality}
T_G(X,Y)\; =\;  T_{G^*}(Y,X)
\end{equation}
where $G$ is a planar graph, and $G^*$ is its dual. (The vertices of $G^*$ correspond to the connected regions
in the complement of $G$ in the plane, and two vertices are connected by an edge in $G^*$ whenever the two corresponding
regions are adjacent.)

In this paper we introduce a $4-$variable polynomial, $P_{G,{\Sigma}}(X,Y,A,B)$, which is an invariant of a graph $G$ embedded
in a closed orientable surface ${\Sigma}$, which satisfies both the contraction-deletion rule and a duality relation analogous to
(\ref{duality}). The variables
$X,Y$ play the same role as in the definition of the Tutte polynomial, while the additional variables $A,B$ reflect
the topological information of $G$ in $\Sigma$. It follows that the Tutte polynomial is a specialization of $P_G$, where
this extra information, reflecting the embedding $G\subset {\Sigma}$, is disregarded.

The main motivation for this work came from an attempt to
understand the combinatorial structure underlying the Potts model on surfaces. As mentioned above, the partition function for the Potts
model on the plane is given by the Tutte polynomial, while on surfaces essential loops are weighted differently from
trivial loops (for references on the Potts model on surfaces, cf. \cite{CM,DSZ,JS,MS}.) This leads to the introduction of additional variables, keeping track of the topological information
of graphs on surfaces.

Using topological duality on surfaces, we establish the duality relation
\begin{equation} \label{duality equation}
P_G\, (X,Y,A,B)\; =\; P_{G^*}\, (Y,X,B,A),
\end{equation}
which may be viewed as a natural analogue of the duality (\ref{duality}) of the Tutte polynomial
for planar graphs.
For the dual graph $G^*$ in (\ref{duality equation}) to be well-defined, it is natural to consider graphs $G$ which are {\em cellulations} of
$\Sigma$, that is, graphs such that each component of ${\Sigma}\smallsetminus G$ is a disk. Equivalently, such graphs
may be viewed as orientable {\em ribbon graphs}, this point of view is presented in more detail in section \ref{ribbon section}.

For ribbon graphs, there is a well-known $3$-variable polynomial defined by B. Bollob\'{a}s and O. Riordan \cite{BR1, BR2}.
We denote this graph polynomial by $BR_G(X,Y,Z)$, its construction is recalled in section \ref{ribbon section}. We show that this polynomial can be
obtained as a specialization of $P_G$:
\begin{equation} \label{relation equation}
BR_G(X,Y,Z)\; =\; Y^{g}\; P_G(X-1,Y,YZ^2,Y^{-1}),
\end{equation}
where $g$ is the genus of the ribbon graph $G$.
In fact the authors prove in \cite{BR1,BR2} that their polynomial is a universal invariant of ribbon graphs with respect to
the contraction-deletion rule (we give a precise statement of this result in section \ref{ribbon section}.) Therefore
in principle the two polynomials $BR_G$, $P_G$ carry equivalent information about the ribbon graph $G$, although
an expression of $P_G$ in terms of $BR_G$ does not seem to be as straightforward as (\ref{relation equation}).
We note that the definition of the polynomial $P_G$ could be normalized so the specialization to $BR_G$ is obtained
by simply setting one of the variables equal to $1$ (see section \ref{ribbon section}). We chose a normalization
making the duality statement (\ref{duality equation}) most natural.

Several authors have established partial results on duality
for the Bollob\'{a}s-Riordan polynomial: B. Bollob\'{a}s and O. Riordan \cite{BR1} stated duality for a $1$-variable
specialization, J. A. Ellis-Monaghan and I. Sarmiento \cite{ES} and I. Moffatt \cite{Mo}
(see also \cite{Ch, ES1, Mo1}) proved duality for a certain $2$-variable specialization.
These results may be recovered as a consequence of equations
(\ref{duality equation}), (\ref{relation equation}), see sections
\ref{duality section} and \ref{prior duality results}; our result (\ref{duality equation}) is more general.

A self-contained discussion of the polynomial $P_{G,{\Sigma}}(X,Y,A,B)$ and its properties may be
found in sections \ref{graph polynomial section},
\ref{duality section 1}. In section \ref{reformulation} we point out a combinatorial formulation of this
 graph polynomial in the context of ribbon graphs, without using homology. The reader interested in a more detailed discussion
of topological aspects of graphs on surfaces
will find in section \ref{linear algebra section} the definition of a more general, infinite-variable, polynomial $\widetilde P_{G,{\Sigma}}$. We point out in that section that a general context for the duality of graph polynomials on an oriented closed surface $\Sigma$
is provided by
the intersection pairing and the Poincar\'{e} duality, giving rise to a symplectic
structure on the first homology group $H_1({\Sigma})$. (The action of the mapping class group of the surface
induces a representation of the symplectic group $Sp\, (2g,{\mathbb Z})$ on $H_1({\Sigma}, {\mathbb Z})$, where $g$ is the genus
of the surface.) Given a subgroup $V$ of $H_1({\Sigma})$, its ``orthogonal complement'' $V^{\perp}$  with respect to
the intersection form may be defined, see (\ref{perp}).
Using this structure, we define in section \ref{linear algebra section} a more general version of the polynomial $P_G$,
with coefficients corresponding to subgroups of $H_1({\Sigma})$, its duality property is stated in Lemma
\ref{duality lemma 2}. This more general polynomial
may be used to distinguish different embeddings of a graph in $\Sigma$. (One may also generalize further and, avoiding
the use of homology, consider the {\em Tutte skein module} of a surface $\Sigma$: the vector space spanned by
isotopy classes of graphs on $\Sigma$, modulo the contraction-deletion relation, see section \ref{Tutte skein}.
In this case the ``polynomial'' associated to a graph $G\subset {\Sigma}$ is the element of the skein
module represented by $G$.)
On the other hand, if one considers
graphs on ${\Sigma}$ up to the action of the diffeomorphism group of $\Sigma$
(or if one studies ribbon graphs), then the relevant invariant is the
finite-variable polynomial $P_G$, discussed above.

In section \ref{Jones section} a version of the Kauffman bracket and of the
Jones polynomial on surfaces is considered, taking into account homological information
on the surface. In particular, using the interpretation of a virtual link
as an ``irreducible'' embedding of a link into a surface due to G. Kuperberg \cite{Ku},
this defines a generalization of the Jones polynomial for
virtual links. For example, the Jones polynomial $J_L(t,Z)$ acquires a new variable $Z$ which,
in the state-sum expression,
keeps track of the rank of the subgroup of the first homology group $H_1({\Sigma})$ of the
surface represented by a resolution of the link diagram on the surface.

If a link $L$ has an alternating
diagram on $\Sigma$, the diagram may be
checkerboard colored, and there is a graph $G$ (the Tait graph) associated to it.
In this context we show (Theorem \ref{Jones Tutte lemma}) that the generalized Kauffman bracket (and the Jones polynomial $J_L(q,Z)$)
is a specialization of the polynomial $P_G$,
generalizing the well-known relation between the Jones polynomial of a link in $3-$space and the
Tutte polynomial associated to its planar projection due to Thistlethwaite \cite{Th}. The analogue
of Thistlethwaite's theorem, relating the Kauffman bracket of virtual links and the Bollob\'{a}s-Riordan polynomial
of ribbon graphs, was established by S. Chmutov and I. Pak in
\cite{ChP}. Theorem \ref{Jones Tutte lemma} generalizes these results to the polynomial
$J_L$ with the extra homological parameter $Z$.
This relation between the generalized Jones polynomial $J_L(q,Z)$ and the polynomial
$P_G$ of the associated graph does not seem to have an immediately obvious analogue in terms of $BR_G$.

A {\em multivariate} version of the Tutte polynomial, where the edges of a graph are weighted, is important
in the analysis of the Potts model \cite{S}. We define its generalization, a multivariate version of the
polynomial $P_G$, and establish a duality analogous to (\ref{duality equation}) in section \ref{multivariable section}.
(A mutlivariate version of the Bollob\'{a}s-Riordan polynomial has been considered by I. Moffatt in \cite{Mo}, and
F. Vignes-Tourneret \cite{VT} established a partial duality
result for a signed version of the multivariate Bollob\'{a}s-Riordan polynomial.)

The Tutte polynomial and the definition of the new polynomial $P_G$, as well
as a discussion of its basic properties, are given in section \ref{graph polynomial section}.
Its duality relation (\ref{duality equation}) is proved in section \ref{duality section 1}.
We review the notion of a ribbon graph and the
definition of the Bollob\'{a}s-Riordan polynomial, and we establish the relation (\ref{relation equation}) in section
\ref{ribbon section}. Section \ref{prior duality results} shows that our duality result (\ref{duality equation})
implies the previously known partial results on duality for the Bollob\'{a}s-Riordan polynomial.
Section \ref{linear algebra section} recalls basic notions of symplectic linear algebra, allowing one to
generalize $P$ to a polynomial $\widetilde P_{G,{\Sigma}}$ with coefficients taking values in subgroups of the first homology group of the
surface.
Section \ref{Jones section} defines the relevant versions of the Jones polynomial and of the Kauffman bracket
and establishes a relationship between them and the polynomial $P_G$.
Finally, section \ref{multivariable section} discusses a multivariate
version of the polynomial $P_G$ and the corresponding duality relation.

{\bf Acknowledgements.}
This work is related to an ongoing project with Paul Fendley \cite{FK1}, \cite{FK2} relating TQFTs, graph polynomials, and algebraic
and combinatorial properties of models of statistical mechanics. I would like to thank Paul for many discussions that motivated the results in this paper.

I would like to thank the referee for the comments on the earlier version of this paper which led to a substantially improved exposition.

\bigskip

\bigskip

\section{The Tutte polynomial and graphs on surfaces} \label{graph polynomial section}

Consider the following normalization of the Tutte polynomial
of a graph $G$:

\begin{equation} \label{Tutte definition}
T_G(X,Y)=\sum_{H\subset G} X^{c(H)-c(G)}\; Y^{n(H)}.
\end{equation}

The summation is taken over all spanning subgraphs $H$ of $G$, that is the subgraphs $H$ such that
the vertex set of $H$ coincides with the
vertex set of $G$. Therefore the sum contains $2^e$ terms, where $e$ is the number of edges of $G$.
In (\ref{Tutte definition}) $c(H)$ denotes the number of connected components of the graph $H$, and
$n(H)$ is the {\em nullity} of $H$, defined as the rank of the first homology group $H_1(H)$ of $H$.
(Note that the nullity $n(H)$ may also be computed as $c(H)+e(H)-v(H)$, where $e$ and $v$ denote the number
of edges and vertices of $H$, respectively.)

Now suppose $G$ is a graph embedded in a surface $\Sigma$. We need to introduce some preliminary topological notions which will
be used in the definition (\ref{poly def equation}) of the graph polynomial below. We note that in the context of ribbon graphs,
there is a formulation of this graph polynomial in purely combinatorial terms, see section \ref{reformulation}.

\begin{definition} \label{s perp}
For a spanning subgraph $H$ of $G$, let $s(H)$ be twice the genus of the surface obtained as a regular neighborhood ${\mathcal H}$ of the graph $H$
in $\Sigma$.
(${\mathcal H}$ is a surface with boundary, and its genus is defined as the genus of the closed surface obtained from ${\mathcal H}$ by
attaching a disk to each boundary circle of ${\mathcal H}$.)
Similarly, let $s^{\perp}(H)$ denote twice the genus of the surface obtained by removing a regular neighborhood ${\mathcal H}$ of $H$ from $\Sigma$.
Denote by $i$ the embedding $G\longrightarrow {\Sigma}$, and consider the induced map on the first homology groups with real coefficients (we mention
\cite{H, Massey} as general references on algebraic topology, in particular for the background on the homology groups). Define
\begin{equation} \label{kernel}
k(H) := {\rm dim}\, ({\rm ker}\, (i_*\co H_1(H;{\mathbb R})\longrightarrow H_1({\Sigma}; {\mathbb R}))).
\end{equation}
\end{definition}
For example, for the graph $H$ on the surface of genus $3$,
consisting of a single vertex and $3$ edges, shown on the left in figure \ref{surface}, $s(H)=s^{\perp}(H)=2$, $k(H)=0$.

Note that $k(H)$, which enters the definition (\ref{poly def equation}) below as the exponent of $Y$,
may be replaced by the nullity $n(H)$ (which is the exponent of $Y$ in the Tutte polynomial (\ref{Tutte definition})), the result would be a
different normalization of the polynomial $P$. See formula
(\ref{useful formula}) in section \ref{ribbon section} below relating $n(H)$ and the parameters used in the definition of $P$.
The choice of the exponent of $Y$ in (\ref{poly def equation}) was motivated by the
duality relation (\ref{duality lemma equation}) which is most naturally stated with this normalization.
We introduce the polynomial $P_{G,{\Sigma}}$
which is the main object of study in this paper:
\begin{equation} \label{poly def equation}
P_{G,{\Sigma}}(X,Y,A,B)\; =\; \sum_{H\subset G} X^{c(H)-c(G)}\; Y^{k(H)}\; A^{s(H)/2}\;
B^{s^{\perp}(H)/2}
\end{equation}
The reader interested in a more general topological context for analyzing polynomial invariants of
graphs on surfaces and their duality properties should compare $P_{G,{\Sigma}}$ with the more general version
defined in section \ref{linear algebra section}. (The invariants $s(H), s^{\perp}(H)$ fit naturally in that context, and
this explains, in part, their normalization as {\em twice} the genus of the corresponding surface.)
Some elementary properties
of the polynomial $P_{G,{\Sigma}}$ are summarized in the following statement.
A surface $\Sigma$ usually will be fixed, and the subscript $\Sigma$ will be omitted from the notation.

\begin{lemma} \label{properties} \sl \nl
(1) If $e$ is an edge of $G$ which is neither a loop nor a bridge, then $P_G=P_{G\smallsetminus e}+P_{G/e}$.\\
(2) If $e$ is a bridge in $G$, then $P_G=(1+X)\, P_{G/e}$.\\
(3) If $e$ is a loop in $G$ which is trivial in $H_1({\Sigma})$,
then $P_G=(1+Y)\, P_{G\smallsetminus e}$.\\
\end{lemma}

{\em Proof.} The proof of lemma \ref{properties} is similar to the proof of the corresponding statements
for the Tutte polynomial. To prove (1), consider an edge $e$ which is neither a loop nor a bridge. Since $e$ is not a loop,
the sum
(\ref{poly def equation}) splits into two parts $P_G=S_1+S_2$. $S_1$ consists of the terms with $H$ containing the edge $e$,
and $S_2$ consists of the terms
with $H$ not containing $e$. In the first case, the embedding $H\subset {\Sigma}$ is homotopic to the
embedding $H/e\subset {\Sigma}$, and all of the invariants $c, k, s, s^{\perp}$ of $H$ coincide with those of $H/e$.
Therefore, $S_1=P_{G/e}$.
The terms in $S_2$ are in $1-1$ correspondence with the terms in $P_{G\smallsetminus e}$. Moreover, since $e$ is not a bridge,
$c(G)=c(G\smallsetminus e)$. It follows that $S_2=P_{G\smallsetminus e}$.

To prove (2), suppose $e$ is a bridge in $G$. Again the sum (\ref{poly def equation}) splits:

$$P_G = \!\!\! \sum_{H\subset (G\smallsetminus e)}\!\!\! X^{c(H)-c(G)}\, Y^{k(H)}\, A^{s(H)/2}\, B^{s^{\perp}(H)/2} \, + \!\! \sum_{H\subset (G/e)}
\!\!\! X^{c(H)-c(G)}\, Y^{k(H)}\, A^{s(H)/2}\, B^{s^{\perp}(H)/2}$$
More precisely, the subgraphs $H$ parametrizing the second sum are all subgraphs of $H$
containing $e$. Contracting $e$ leaves each term in the second sum unchanged, and moreover the second sum
is precisely the expansion of $P_{G/e}$.

There is a $1-1$ correspondence between the subgraphs $H$ (not containing $e$) of $G$ parametrizing the first sum and the subgraphs
parametrizing the second sum. Given $H\subset G$, $e\notin H$, this correspondence associates to
it the subgraph $\widetilde H\subset G/e$ obtained by identifying the two endpoints of $e$ in $H$.
Since $e$ is a bridge, the homological invariants $k, s, s^{\perp}$ of $H$ are identical to those
of $\widetilde H$. However $c(H)-c(G)=c(\widetilde H)-c(G/e)+1$. Therefore each term in the first sum
equals the corresponding term in the expansion of $P_{G/e}$ times $X$. This concludes the proof of (2).

The proof of (3) is analogous, noting that removing a loop $e$ which is homologically trivial on the
surface from a subgraph $H$ decreases $k(H)$ by $1$ and leaves other exponents in the expansion
(\ref{poly def equation}) unchanged. The proof that the exponent of $B$ does not change relies on the fact
that deleting a homologically trivial loop $a$ disconnects the surface. This fact conceptually is a consequence of the
Poincar\'{e}-Lefschetz duality \cite{H}. It may also be observed using a more elementary argument as follows. Supposing the opposite is true, one
immediately finds a loop $b$ in ${\Sigma}$ which intersects $a$ in a single point. This is a contradiction
with the basic and fundamental fact in homology theory that a boundary has trivial intersection number
with any cycle. (See \cite{Massey} for a general discussion of curves on surfaces as well as for applications of homology
in this context.)
\qed

\begin{figure}[ht]
\medskip
\includegraphics[height=2.8cm]{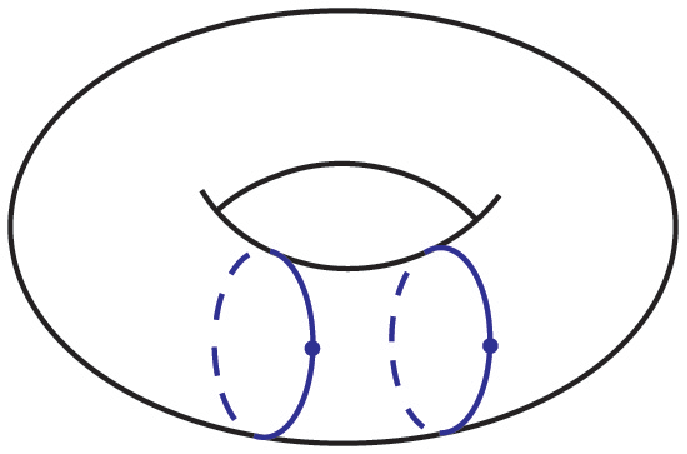}
{\scriptsize
    \put(-98,20){$G_1$}
    \put(-30,20){$G_2$}}
\medskip
\caption{In general the polynomial $P_G$ is not multiplicative with respect to
disjoint unions: in this example, $$P_{G_1}=P_{G_2}=1+B,\;  P_{G_1\cup G_2}=2+B+Y.$$ Note
that $P_G$ is multiplicative for {\em ribbon} graphs, see section \ref{ribbon section}.}
\label{torus}
\end{figure}

{\bf Remark.} Note that if $G_1$, $G_2$ are disjoint graphs in $\Sigma$, it is {\em not} true
in general that  $P_{{G_1}\cup {G_2}}=P_{G_1}\, P_{G_2}$, see for example figure \ref{torus}.
(A similar comment applies to the case when $G_1,G_2$ in $\Sigma$ are disjoint except for a single vertex $v$.)
This is quite different form the case of the classical Tutte polynomial. However, the polynomial $P_G$ is multiplicative
with respect to disjoint unions in the context of ribbon graphs, see lemma \ref{disjoint ribbon}.
The proof of that lemma
also shows that if $G_1$, $G_2$ are graphs in disjoint surfaces, $G_1\subset {\Sigma}_1$, $G_2\subset {\Sigma}_2$,
then $P_{G_1\sqcup G_2, {\Sigma}_1\sqcup {\Sigma}_2}=P_{G_1, {\Sigma}_1}\,
P_{G_2, {\Sigma}_2}$.

\begin{lemma} \label{relation with Tutte} \sl The Tutte polynomial (\ref{Tutte definition}) is a specialization of $P_G$:
$$ T_G(X,Y)\; = \; Y^{g}\; P_{G,{\Sigma}}(X,Y,Y,Y^{-1}),$$
where $g$ is the genus of the surface $\Sigma$.
\end{lemma}

{\em Proof.}
Substituting $A=Y$, $B=Y^{-1}$ into the expansion
(\ref{poly def equation}), one
gets terms of the form
\begin{equation}\label{poly terms}
X^{c(H)-c(G)}Y^{k(H)+s(H)/2-s^{\perp}(H)/2}.
\end{equation}
We claim that
\begin{equation} \label{eq: useful identity}
n(H)=k(H)+g+s(H)/2-s^{\perp}(H)/2.
\end{equation}
This formula shows that each term of the form (\ref{poly terms}) above, multiplied by
$Y^g$, gives the corresponding term in the expansion (\ref{Tutte definition}) of the Tutte polynomial.

The conceptual framework for the formula (\ref{eq: useful identity}) is provided by the structure
of the first homology group $H_1({\Sigma})$ of the surface given by the intersection numbers of curves in
$\Sigma$. This point of view is discussed
in more detail in section \ref{linear algebra section} of the paper, specifically see the identities
(\ref{identities}). At outline of the
argument may be seen as follows.
Following definition \ref{s perp}, given a spanning subgraph $H$ of $G$ let ${\mathcal H}$ be its regular neighborhood in $\Sigma$.
Denote by ${\mathcal H}^*$ its complement: ${\mathcal H}^*={\Sigma}\smallsetminus {\mathcal H}$, so $\Sigma$ is represented as
the union of two compact surfaces ${\mathcal H},{\mathcal H}^*$ along their boundary.
The dimension $2g$ of $H_1({\Sigma})$ equals $s(H)+s^{\perp}(H)+2l$ for some integer $l\geq 0$.
Observe that the dimension of the image of $H_1({\mathcal H})$ in $H_1({\Sigma})$ is precisely $s(H)+l$. 
Indeed, if the dimension
of the image were greater than $s(H)+l$ then in fact the genus of ${\mathcal H}$ must be greater than $s(H)/2$, a contradiction.
Similarly,
if the dimension of the image were less than $s(H)+l$ then the genus of ${\mathcal H}^*$ must be greater than $s^{\perp}(H)/2$.
This shows that $2g=s(H)+s^{\perp}(H)+2l$, and together with $n(H)=k(H)+s(H)+l$ this establishes
(\ref{eq: useful identity}), concluding the proof of lemma \ref{relation with Tutte}. 
See section \ref{linear algebra section} for a more detailed discussion of the
underlying homological structure. We note that a less direct combinatorial proof of the formula (\ref{eq: useful identity}) may be given
using the combinatorial interpretation of the invariants $s(H), s^{\perp}(H)$. (In particular, $s(H)=
c(H)-bc(H)+n(H)$, see section \ref{reformulation}.)
 \qed

{\bf Remark.} The polynomial $P_G$ can be normalized to make the relation with the Tutte polynomial easier to state.
For example, if one chose the exponent of $Y$ in (\ref{poly def equation}) to be $n(H)={\rm rank}\, H_1(H)$ rather than
$k(H)$, $T_G(X,Y)$ would be the specialization of the resulting polynomial obtained simply by setting $A=B=1$.
We chose the convention ({\ref{poly def equation}) to have a natural expression of duality (\ref{duality equation}),
proved in theorem \ref{duality lemma}
below.

\bigskip

\section{Duality.} \label{duality section 1}
In this section we prove a duality result for the polynomial $P_G$ defined in (\ref{poly def equation}),
which is analogous to the duality $T_G(X,Y)=T_{G^*}(Y,X)$ satisfied by the Tutte polynomial
of planar graphs. The following result applies to {\em cellulations} of surfaces: graphs $G\subset {\Sigma}$ such that
each connected component of ${\Sigma}\smallsetminus G$ is a disk. This is a natural condition guaranteeing that
the dual $G^*$ is well-defined. Equivalently, one may view $G$ as a {\em ribbon graph}, see section \ref{ribbon section}.

\begin{theorem} \label{duality lemma} \sl
Suppose $G$ is a cellulation of a closed orientable surface $\Sigma$ (equivalently, let $G$ be an oriented
ribbon graph.)
Then the polynomial invariants of $G$ and its dual $G^*$ are related by
\begin{equation} \label{duality lemma equation}
P_G\, (X,Y,A,B)\; =\; P_{G^*}\, (Y,X,B,A)
\end{equation}
\end{theorem}

{\em Proof.} Consider the expansions (\ref{poly def equation}) of both sides in the statement of the theorem:
\begin{equation} \label{sum1}
P_{G}(X,Y,A,B)\; =\; \sum_{H\subset G} X^{c(H)-c(G)}\; Y^{k(H)}\; A^{s(H)/2}\; B^{s^{\perp}(H)/2}
\end{equation}
\begin{equation}\label{sum2}
P_{G^*}(Y,X,B,A)\; =\; \sum_{H^*\subset G^*} Y^{c(H^*)-c(G^*)}\; X^{k(H^*)}\; B^{s(H^*)/2}\; A^{s^{\perp}(H^*)/2}.
\end{equation}
Recall that the vertices of $G^*$ correspond to the components of the complement ${\Sigma}\smallsetminus G$ (which
are disks since $G$ is a cellulation), and two vertices are connected by an edge in $G^*$ if and only if the
two corresponding components share an edge. Therefore the edges of $G$ and $G^*$ are in $1-1$ correspondence,
with each edge $e$ of $G$ intersecting the corresponding edge $e^*$ of $G^*$ in a single point, and $e$ is disjoint
from all other edges of $G^*$.
For each spanning subgraph $H\subset G$, consider the spanning subgraph $H^*\subset G^*$ whose edges are precisely all
those edges of $G^*$ which are disjoint from all edges of $H$. The theorem follows from the claim that the term corresponding
to $H$ in the expansion of $P_G$ equals the term corresponding to $H^*$ in the expansion of $P_{G^*}$
above.

\begin{figure}[t]
\medskip
\includegraphics[height=3cm]{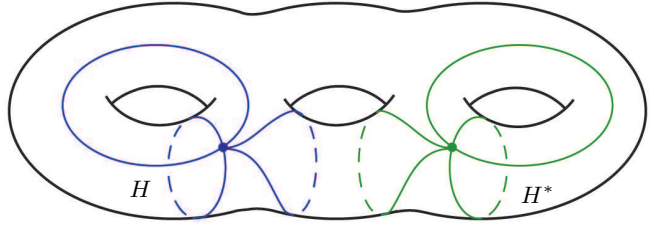}
{\scriptsize
    \put(-198,10){$H$}
    \put(-50,8){$H^*$}}
\medskip
\caption{A subgraph $H$ of a cellulation $G$ of the genus $3$ surface $\Sigma$, and the corresponding
subgraph $H^*$ of the dual cellulation $G^*$.}
\label{surface}
\end{figure}

The simplest example of a cellulation $G$ of a surface $\Sigma$ of genus $g$ is a graph consisting of
a single vertex and $2g$ edges which are loops representing a symplectic basis of $H_1({\Sigma})$. Then its
dual is the graph $G^*$ also with a single vertex and $2g$ loops.
Figure \ref{surface}
shows the surface of genus $3$ and a subgraph $H\subset G$ formed by $3$ edges on the left in the figure. In this case $H^*$ also
consists of $3$ edges as
illustrated on the right in the figure.

The two cellulations $G, G^*$ give rise to
dual handle decompositions of the surface ${\Sigma}$. In the handle decomposition corresponding to $G$, the $0-$handles
are disk neighborhoods of the vertices of $G$, the $1-$handles are regular neighborhoods of the edges of
$G$, the $2-$handles correspond to the $2-$cells ${\Sigma}\smallsetminus G$. Let ${\mathcal H}$ (respectively
${\mathcal H}^*$) denote the union of the $0-$ and $1-$handles corresponding to $H$ (respectively $H^*$).
Note that ${\mathcal H}$ is a regular neighborhood of the graph $H$, and similarly ${\mathcal H}^*$ is
a regular neighborhood of the graph $H^*$.
If $H$ is the entire graph $G$, $H^*$
consists of all vertices of $G^*$ and no edges.
Removing one edge from $H$ at a time, observe that the effect on ${\mathcal H}$ is the removal of a
$1-$handle, while the effect on the dual handle decomposition ${\mathcal H}^*$ is the addition of the co-core of the removed $1-$handle.
To summarize, $\Sigma$ is the union of two surfaces $\mathcal H$, ${\mathcal H}^*$ along their boundary, where ${\mathcal H}$ is
a regular neighborhood of $H$, and ${\mathcal H}^*$ is a regular neighborhood of $H^*$.

It follows from definition of $s, s^{\perp}$ that $s(H)=s^{\perp}(H^*)$ and $s^{\perp}(H)=
s(H^*)$. One also checks that $c(H^*)-c(G^*)=k(H)$ and $c(H)-c(G)=k(H^*)$. A geometric argument may be given for this fact,
where one considers the induction on the number of edges in $H$ and observes that adding an edge to $H$ (and therefore
removing an edge from $H^*$) either
decreases both $c(H)-c(G)$ and $k(H^*)$ by one, or leaves both quantities unchanged. We give a more direct, algebraic-topological
proof: by Poincar\'{e}-Lefschetz duality (cf \cite{H}), since ${\Sigma}={\mathcal H}\cup {\mathcal H}^*$,
one has an isomorphism of the relative second homology
group $H_2({\Sigma}, {\mathcal H})$ and the $0$th cohomology $H^0({\mathcal H}^*)$. The dimension of
$H^0({\mathcal H}^*)\cong H_0({\mathcal H}^*)$ equals $c(H^*)$, the number of connected components of $H^*$.
The group $H_2({\Sigma}, {\mathcal H})$ fits in the long exact sequence (cf \cite{H}) of the pair $({\Sigma}, {\mathcal H})$:
$$ 0\longrightarrow H_2({\Sigma})\longrightarrow H_2({\Sigma}, {\mathcal H})\longrightarrow H_1({\mathcal H})
\longrightarrow H_1({\Sigma})\longrightarrow \ldots,$$
therefore $$dim(H_2({\Sigma}, {\mathcal H}))=dim(ker[H_1({\mathcal H})\longrightarrow H_1({\Sigma})])+dim(H_2({\Sigma})).$$
The first term in this sum is the definition of $k(H)$, and each connected component of ${\Sigma}$ contributes $1$ to
the second term. Since $G^*$ is a cellulation, the number of connected components of $G^*$ equals the number of connected components
of $\Sigma$. Combining these equalities, one has $c(H^*)=k(H)+c(G^*)$. This proves $c(H^*)-c(G^*)=k(H)$, and analogously
one has $c(H)-c(G)=k(H^*)$.

This shows that the
terms corresponding to $H, H^*$ in (\ref{sum1}, \ref{sum2}) are equal, concluding the proof of theorem \ref{duality lemma}.
\qed

Note that duality results for certain specializations of the Bollob\'{a}s-Riordan polynomial have been previously obtained by several authors.
We discuss these results and show that they may be derived as a consequence of our theorem \ref{duality lemma} in section
\ref{prior duality results} below.

\bigskip

\section{Ribbon graphs and the Bollob\'{a}s-Riordan polynomial} \label{ribbon section}

A {\em ribbon graph} is a pair $(G,S)$ where $G$ is a graph embedded in a surface (with boundary) $S$ such that
the embedding $G\hookrightarrow S$ is a homotopy equivalence. It is convenient
to consider the surface $S$ with a handle decomposition corresponding to the graph $G$: the $0$-handles are disk
neighborhoods of the vertices of $G$, and the $1$-handles correspond to regular neighborhoods of the edges.
(Other terms: cyclic graphs, fat graphs are also sometimes used
in the literature to describe ribbon graphs.)
$G$ is an {\em orientable} ribbon graph if
$S$ is an orientable surface. Given a ribbon graph $(G,S)$, one obtains a closed surface $\Sigma$ by attaching
a disk to $S$ along each boundary component.
Therefore a ribbon graph may be viewed as a {\em cellulation} of a closed surface
$\Sigma$, i.e. a graph $G$ embedded in $\Sigma$ such that each component of the complement
${\Sigma}\smallsetminus G$ is a disk. Conversely, given a cellulation $G$ of $\Sigma$, one has a ribbon graph structure
$(G,S)$ where $S$ is a regular neighborhood of $G$ in $\Sigma$.
 We will use the notions
of a ribbon graph and of a cellulation interchangeably.

Consider the Bollob\'{a}s-Riordan polynomial of ribbon graphs \cite{BR1, BR2} (in this paper
we only consider {\em orientable} ribbon graphs, therefore there are three, rather than four, variables): given a ribbon
graph $(G,S)$,

\begin{equation} \label{BR}
BR_{G,S}(X,Y,Z)\; =\; \sum_{H\subset G} (X-1)^{r(G)-r(H)}\, Y^{n(H)}\, Z^{c(H)-bc(H)+n(H)}.
\end{equation}

The summation is taken over all spanning subgraphs $H$ of $G$, and moreover each $H$ inherits the ribbon structure from that of $G$:
the relevant surface is obtained as the union of all $0$-handles and just those $1$-handles which correspond to the edges of $H$.
To explain the notation in this definition, let $v(H), e(H)$ denote the number of vertices, respectively edges, of $H$,
and let $c(H)$ be the number of connected components. ($v(H)=v(G)$ since $H$ is a spanning subgraph of $G$.)
Then $r(H)=v(G)-c(H)$, $n(H)=e(H)-r(H)$, and $bc(H)$ is the number of boundary components of the surface $S$. Note that
$n(H)$ equals the rank of the first homology group $H_1(H)$, and the exponent of $Z$, $c(H)-bc(H)+n(H)$, equals $2g(H)=s(H)$, twice
the genus of the surface underlying the ribbon graph $H$. To simplify the notation, we will often omit the reference to the surface $S$
and denote the polynomial by $BR_G$.

\begin{lemma} \label{ribbon lemma} \sl
The Bollob\'{a}s-Riordan polynomial of a ribbon graph may be obtained as a specialization
of the polynomial $P_G$:
\begin{equation}
BR_{G,S}(X,Y,Z)\; =\; Y^{g}\; P_{G,{\Sigma}}(X-1,Y,YZ^2,Y^{-1}),
\end{equation}
where ${\Sigma}$ is the closed surface obtained by attaching
a disk to $S$ along each boundary component, and $g$ is the genus of $\Sigma$.
\end{lemma}

The proof consists of showing that the corresponding terms in the expansions (\ref{poly def equation}), (\ref{BR}) are equal.
Indeed, substituting  $A=YZ^2$, $B=Y^{-1}$ in (\ref{poly def equation}) gives summands of the form
$$(X-1)^{c(H)-c(G)}\; Y^{k(H)+s(H)/2-s^{\perp}(H)/2}\; Z^{s(H)}.$$
Using the formulas (\ref{identities}), established in the following section, observe that
$$n(H)=k(H)+g+s(H)/2-s^{\perp}(H)/2,$$
therefore these summands are equal to $Y^{-g}$ times the corresponding terms in
(\ref{BR}).
\qed

To discuss the relation between the polynomial $P_G$ and the Bollob\'{a}s-Riordan polynomial further, recall that the polynomial $BR_G$ satisfies the following universality property. Let ${\mathcal G}$ denote the set of isomorphism classes \cite{BR2} of connected ribbon graphs.
Define the maps $C_{ij}$ from ${\mathcal G}$ to ${\mathbb Z}[X]$ by $BR=\sum_{i,j} C_{ij} Y^iZ^j$. Further, given a commutative
ring $R$ and an element $x\in R$, $C_{ij}(x)$ will
denote the map from ${\mathcal G}$ to $R$ obtained by composing $C_{ij}$ with the ring homomorphism
${\mathbb Z}[X]\longrightarrow R$ mapping $X$ to $x$.

\begin{theorem} \label{universality} \cite{BR1,BR2} \sl
Let $R$ be a commutative ring and $x\in R$ and  ${\phi}\co {\mathcal G}\longrightarrow R$ a map satisfying\\
(1) ${\phi}(G)={\phi}(G/e)+{\phi}(G\smallsetminus e)$ if $e$ is neither a loop nor a bridge, and \\
(2) ${\phi}(G)=x\, {\phi}(G/e)$ if $e$ is a bridge.

Then there are elements ${\lambda}_{ij}\in R$, $0\leq j\leq i$, such that
$${\phi}\; =\; \sum_{i,j} {\lambda}_{ij} C_{ij}(x).$$
\end{theorem}

The polynomial $P_G$ satisfies the properties (1), (2) in this theorem, therefore it follows
that the coefficients of $P_G$ may be expressed as linear combinations of the coefficients of $BR_G$.
The main
difference in the definitions of the two polynomials is that each term in the expansion (\ref{BR}) of $BR_G$
is defined in terms of the invariants of a ribbon subgraph $H$, while the terms in the expansion
(\ref{poly def equation}) involve the invariants associated to the embedding of $H$ into the original fixed
surface $\Sigma$.
Indeed, note that the four parameters $c(h), k(H), s(H), s^{\perp}(H)$ in the
definition (\ref{poly def equation}) of $P_G$ are independent invariants of $H$, in the sense that there
are examples of graphs showing that no three of the parameters determine the other one.
Therefore it does not seem likely
that there is a straightforward expression for $P_G$ in terms of $BR_G$ similar to that in lemma \ref{ribbon lemma},
however it would be interesting to find an explicit expression.

Returning to the properties of the polynomial $P_G$, observe that the multiplicativity for disjoint unions and
for one-point unions holds in the context of ribbon graphs (compare with the remark after lemma \ref{properties}):

\begin{lemma} \label{disjoint ribbon} \sl
Properties (1)--(3) in lemma \ref{properties} hold for ribbon graphs $G$. In addition, for disjoint ribbon graphs
$G_1, G_2$,

(4) $P_{{G_1}\sqcup {G_2}}\; =\; P_{{G_1}\vee {G_2}}\; =\; P_G\cdot P_{G'}$.
\end{lemma}

Here by the polynomial $P_G$ of a ribbon graph $(G,S)$ we mean $P_{G,{\Sigma}}$ where as above $\Sigma$ is the closed surface
associated to $S$. For example, the closed surface associated to the graphs $G_1, G_2$ with the ribbon structure
inherited from their embedding into the torus in figure \ref{torus} is the $2$-sphere
(and the surface associated to $G_1\cup G_2$ is the disjoint union of two spheres), and not the torus. This illustrates the difference
between the validity of the property (4) for ribbon graphs, but not in general for graphs on surfaces as in figure \ref{torus}.

{\em Proof.} The proof of (1)--(3) is identical to that in lemma \ref{properties}.
To prove (4) for the disjoint union $G_1\sqcup G_2$, consider subgraphs $H_1\subset G_1$, $H_2\subset G_2$ and
let $V_i$ denote the image of $H_1(H_i)$ in
$H_1({\Sigma})$, $i=1,2$. Since the surface
associated to $G_1\sqcup G_2$ is the disjoint union of surfaces associated to $G_1$ and $G_2$,
one has $k(H_1\cup H_2)=k(H_1)+k(H_2), s(H_1\cup H_2)=s(H_1)+s(H_2)$, and $s^{\perp}(H_1\cup H_2)=s^{\perp}(H_1)+s^{\perp}(H_2)$.
The proof for the one-vertex union $G_1\vee G_2$ is directly analogous.
\qed

\subsection{Prior results on duality of the Bollob\'{a}s-Riordan polynomial} \label{prior duality results}
Several authors have established duality for certain specializations of the Bollob\'{a}s-Riordan polynomial.
For example, \cite{BR2} notes that

\begin{equation} \label{partial duality1}
BR_G(1+t, t, t^{-1})\; =\; BR_{G^*}(1+t, t, t^{-1}).
\end{equation}

By lemma \ref{ribbon lemma}, $BR_G(1+t, t, t^{-1})=Y^g\, P_G(t,t,t^{-1},t^{-1})$, therefore (\ref{partial duality1}) is a
consequence of theorem \ref{duality lemma}. More generally, it is shown in \cite{ES,Mo} (see also \cite{Ch, ES1, Mo1}) that there is duality
for a $2-$variable specialization:

\begin{equation} \label{partial duality2}
BR_G(1+X,Y,(XY)^{-1/2})\; =\; (X^{-1}Y)^g\; BR_{G^*}(1+Y,X,(XY)^{-1/2})
\end{equation}

Observe that according to lemma
\ref{ribbon lemma},

$$ BR_G(1+X,Y,(XY)^{-1/2})\; =\; Y^g\; P_G(X,Y,X^{-1}, Y^{-1}),$$
$$ BR_{G^*}(1+Y,X,(XY)^{-1/2})\; =\;
X^g\; P_{G^*}(Y,X,Y^{-1}, X^{-1}).$$
Therefore (\ref{partial duality2}) may also be viewed as a consequence of theorem \ref{duality lemma}. It would be
interesting to understand the full duality relation (\ref{duality lemma equation}) in terms of the Bollob\'{a}s-Riordan polynomial,
since as discussed above, the polynomials $P_G, BR_G$ carry equivalent information about a ribbon graph $G$.

\subsection{A reformulation of the polynomial ${\mathbf{P_G}}$ for ribbon graphs.} \label{reformulation}
We conclude this section by noting that one may give a combinatorial formulation of
the polynomial $P_{G}$ (defined by (\ref{poly def equation})) in the context of ribbon graphs.

Given a ribbon graph $(G,S)$, as above consider $G$ as a cellulation of a closed surface $\Sigma$, so each
component of ${\Sigma}\smallsetminus G$ is a disk.
The dual cellulation $G^*\subset{\Sigma}$ is then well-defined: the vertices of $G^*$ correspond to the
components of the complement ${\Sigma}\smallsetminus G$, and two vertices are connected by an edge in $G^*$ if and only if
the corresponding components of ${\Sigma}\smallsetminus G$ share an edge. Therefore the edges of $G$ and $G^*$ are in $1-1$ correspondence,
with each edge $e$ of $G$ intersecting the corresponding edge $e^*$ of $G^*$ in a single point, and $e$ is disjoint
from all other edges of $G^*$.
For each spanning subgraph $H\subset G$, consider the spanning subgraph $H^*\subset G^*$ whose edges are precisely all
those edges of $G^*$ which are disjoint from all edges of $H$.

Given a  ribbon graph $(G,S)$, consider
\begin{equation} \label{preliminary}
P'_{G,S}(X,Y,A,B)\; =\; \sum_{H\subset G} X^{c(H)-c(G)}\, Y^{n(H)}\, A^{c(H)-bc(H)+n(H)}B^{c(H^*)-bc(H^*)+n(H^*)}.
\end{equation}
The summation is taken over all spanning ribbon subgraphs $H$ of $G$. Note that
the exponent of $A$, $c(H)-bc(H)+n(H)$, equals $2g(H)$, twice
the genus of the surface underlying the ribbon graph $H$. Similarly, the exponent of $B$ equals twice
the genus of the dual ribbon graph $H^*$. Since these quantities correspond to the invariants $s(H), s^{\perp}(H)$ (see definition \ref{s perp}), the polynomial $P'$ may be rewritten as
\begin{equation}
P'_{G,S}(X,Y,A,B)\; =\; \sum_{H\subset G} X^{c(H)-c(G)}\, Y^{n(H)}\, A^{s(H)}B^{s^{\perp}(H)}.
\end{equation}
The invariants $n(H), k(H), s(H), s^{\perp}(H)$ may be related using the formulas (\ref{identities}), established in the following section:
\begin{equation} \label{useful formula}
n(H)=k(H)+g+s(H)/2-s^{\perp}(H)/2,
\end{equation}
where $g$ is the genus of the surface underlying the ribbon graph $G$.
Then it is straightforward to see that the polynomial $P_G$, defined by (\ref{poly def equation}), and $P'_G$, combinatorially defined above, are equivalent (may be obtained from each other by substitutions of variables). For example, $P_G$ may be expressed in terms of $P'_G$ as follows:
$$P_G(X,Y,A,B)=Y^{-g}\, P'_G(X,Y,A Y^{-1/2}, B Y^{1/2}).$$

\bigskip

\section{Symplectic linear algebra and a more general graph polynomial} \label{linear algebra section}

In this section we show that the polynomial $P_{G,{\Sigma}}$, defined in section \ref{graph polynomial section},
fits in a more general topological framework. (The material in this section is not directly used in sections \ref{Jones section},
\ref{multivariable section}, and therefore the reader who is interested in applications to knot theory, or in the multivariate
version of the graph polynomial $P_G$, may  choose to proceed directly to the subsequent sections of the paper.)
First we recall a number of basic facts and introduce certain notation in the
symplectic linear algebra setting which will be useful for the definition of
the more general graph polynomial $\widetilde P_{G,{\Sigma}}$. Let $\Sigma$ be
a (not necessarily connected) closed oriented surface, and consider the intersection pairing
$$w\co H_1({\Sigma}, {\mathbb R})\times H_1({\Sigma},{\mathbb R})\longrightarrow {\mathbb R}.$$
The intersection pairing may be viewed geometrically, as the intersection number (where the intersection points are counted with signs) of oriented cycles representing homology classes in $H_1({\Sigma})$, or dually as the cup product on first cohomology $H_1({\Sigma})$, see \cite{H}.
The invariants
considered below do not depend on the orientation.
Poincar\'{e} duality \cite{H} implies that the bilinear form $w$ is non-degenerate, in other words it is
a symplectic form on the vector space $H_1({\Sigma}, {\mathbb R})$. A note on the homology coefficients:
the invariants below may be defined using either ${\mathbb Z}$ or ${\mathbb R}$, and these coefficients will be used interchangeably.

Let $H$ be a graph embedded in the surface $\Sigma$,
and let $i\co H\hookrightarrow {\Sigma}$ denote the embedding.
Denote
\begin{equation} \label{V}
V=V(H)={\rm image}\, (i_*\co H_1(H;{\mathbb R})\longrightarrow H_1({\Sigma}; {\mathbb R})).
\end{equation}
In other words, $V$ is the subgroup of the first homology group of the surface, generated by the cycles in the graph $G$.
The ``symplectic orthogonal complement'' of $V$ may be defined by
\begin{equation} \label{perp}
V^{\perp}=V^{\perp}(H)\; =\; \{ u\in H_1({\Sigma}, {\mathbb R})|\, \forall v\in V(H), \, w(u,v)=0\}.
\end{equation}
The invariants $s(H),\,  s^{\perp}(H)$ of a graph $H$ on ${\Sigma}$, introduced in definition \ref{s perp}, may be defined in this framework as follows:
\begin{equation} \label{symplectic invariants}
s(H)= {\rm dim}(V/(V\cap V^{\perp})), \; \; s^{\perp}(H) = {\rm dim}(V^{\perp}/(V\cap V^{\perp})).
\end{equation}
Said differently, $s(H)$ is the dimension of a maximal symplectic subspace of $V$ (with respect to the
symplectic form $w$ on $H_1({\Sigma}, {\mathbb R})$), and similarly $s^{\perp}(H)$ is the dimension of a maximal symplectic subspace
in $V^{\perp}$. (The fact that (\ref{symplectic invariants}) gives the same invariants as definition \ref{s perp} may be observed
by considering a regular neighborhood ${\mathcal H}$ of $H$ in $\Sigma$ and noting that the homology classes corresponding to
the boundary curves of ${\mathcal H}$ in $H_1({\Sigma})$ are in the intersection $V\cap V^{\perp}$.)
Also it will be useful to consider
\begin{equation} \label{lagrange invariants}
l(H) := {\rm dim}\, (V\cap V^{\perp}), \;\;
k(H) := {\rm dim}\, ({\rm ker}\, (i_*\co H_1(H;{\mathbb R})\longrightarrow H_1({\Sigma}; {\mathbb R}))).
\end{equation}
Note the identities relating these invariants for any graph $H\subset {\Sigma}$:
\begin{equation} \label{identities}
s(H)+s^{\perp}(H)+2l(H)\; =\; 2g, \; \; k(H)+l(H)+s(H) \; =\; {\rm dim}\, (H_1(H)),
\end{equation}
where $g$ denotes the genus of $\Sigma$.

\subsection{A more general graph polynomial} \label{general polynomial}

Now suppose $G$ is a graph embedded in a surface $\Sigma$, let $i\co G\longrightarrow {\Sigma}$
denote the embedding. Consider a collection of formal variables corresponding to the subgroups
of $H_1({\Sigma})$. Given a subgroup $V< H_1({\Sigma})$, let $[V]$ denote the
corresponding variable associated to it. Define

\begin{equation} \label{multi poly def equation}
\widetilde P_{G,{\Sigma}}(X,Y)\; =\; \sum_{H\subset G} [i_*( H_1(H))]\;
X^{c(H)-c(G)}\; Y^{k(H)}.
\end{equation}

Here $[i_*( H_1(H))]$ is the formal variable associated to the subgroup equal to the image
of $H_1(H)$ in $H_1({\Sigma})$ under the homomorphism $i_*$ induced by inclusion; $k(H)$ is defined in
(\ref{kernel}). Therefore $\widetilde P_{G,{\Sigma}}$ may be viewed as a polynomial in $X,Y$ with
coefficients corresponding to the subgroups of $H_1({\Sigma})$. This polynomial may be used to distinguish
different embeddings of a graph $G$ into $\Sigma$.

However if two graphs $G,G'$ in ${\Sigma}$ are considered equivalent whenever there
is a diffeomorphism taking $G$ to $G'$, one needs to consider a polynomial invariant in terms of
quantities which are invariant under the action of the mapping class group. This is the context in which the polynomial
$P_{G,{\Sigma}}$ (defined in section \ref{graph polynomial section}) is useful, indeed it may be viewed as a specialization of $\widetilde P_{G,{\Sigma}}$ where
$[i_*( H_1(H))]$ is specialized to $A^{s(H)/2}B^{s^{\perp}(H)/2}$. In section \ref{duality section} below we point out that
the polynomial $\widetilde P_{G,{\Sigma}}$ satisfies a natural duality relation, generalizing theorem \ref{duality lemma}.

\subsection{The Tutte skein module}  \label{Tutte skein}
One may generalize the polynomial $P_G$ further and, avoiding
the use of homology, consider the {\em Tutte skein module} of a surface $\Sigma$: the vector space spanned by
isotopy classes of graphs on $\Sigma$, modulo relations (1)-(3) in lemma \ref{properties}. For example, the contraction-deletion relation
states that $G=G\smallsetminus e+G/e$, where the three graphs $G, G\smallsetminus e, G/e$ are viewed as vectors in the skein module.
In this case the ``polynomial'' associated to a graph $G\subset {\Sigma}$ is the element of the skein
module represented by $G$. There is an expansion, analogous to (\ref{multi poly def equation}), where each term in the expansion
is an element of the skein module, and to get the polynomial $\widetilde P_{G,{\Sigma}}$ one applies homology to that expansion.

Note that a relative version of this skein module, specialized to $Y=0$, in the rectangle -- the {\em chromatic algebra} --  was considered in
\cite{FK1, FK2}. See also remark 6 following the statement of theorem
\ref{Jones Tutte lemma} below concerning the relation between the Tutte skein module of $\Sigma$ and the Kauffman skein module of
${\Sigma}\times I$.

\subsection{Duality} \label{duality section}

In the remaining part of this section we show that
the polynomial $\widetilde P_{G,{\Sigma}}$ satisfies a natural duality relation, generalizing theorem \ref{duality lemma}:

\begin{lemma} \label{duality lemma 2} \sl
Suppose $G$ is a cellulation of a closed orientable surface $\Sigma$ (equivalently, let $G$ be an oriented
ribbon graph.)
Then $\widetilde P_{G^*}(Y,X)$ is obtained from $\widetilde P_G(X,Y)$ by
replacing each coefficient $[V]$ (formally corresponding to a subgroup of $H_1({\Sigma})$) with its symplectic
orthogonal complement $[V^{\perp}]$.
\end{lemma}

The proof of this lemma follows along the lines of the proof of theorem \ref{duality lemma},
one shows that each term $[i_*( H_1(H))]\; X^{c(H)-c(G)}\; Y^{k(H)}$ in the expansion of $\widetilde P_G(X,Y)$
equals the corresponding term $[i_*( H_1(H^*))]\; X^{k(H^*)}\; Y^{c(H^*)-c(G^*))}$ in the expansion of $\widetilde P_{G^*}(Y,X)$,
and moreover that $i_*( H_1(H^*))\cong (i_*( H_1(H)))^{\perp}$. Here for each spanning subgraph $H\subset G$,
$H^*$ is the ``dual'' subgraph of $G^*$.

The proof of theorem \ref{duality lemma} established that $c(H^*)-c(G^*)=k(H)$ and $c(H)-c(G)=k(H^*)$. The remaining step is
to show that, in the notation of (\ref{V}), (\ref{perp}),
\begin{equation} \label{observation}
V(H^*)\cong V(H)^{\perp}.
\end{equation}

Consider the regular neighborhoods ${\mathcal H}, {\mathcal H}^*$ of $H, H^*$ in $\Sigma$. Then $\Sigma$ is represented as the union
of two surfaces ${\mathcal H}, {\mathcal H}^*$ along their boundary. Since the intersection of any $1$-cycle
in ${\mathcal H}$ with any $1$-cycle in ${\mathcal H}^*$ is zero, it is clear that $V(H^*)\subset V(H)^{\perp}$.
To prove the opposite inclusion, consider the Mayer-Vietoris sequence \cite{H}:
$$\ldots\longrightarrow H_1({\mathcal H})\oplus H_1({\mathcal H}^*)\overset{\alpha}{\longrightarrow} H_1({\Sigma})\overset{\beta}{\longrightarrow} H_0(\partial)\longrightarrow \ldots$$
where $\partial$ denotes $\partial {\mathcal H}=\partial {\mathcal H}^*={\mathcal H}\cap {\mathcal H}^*$. It follows from the geometric
 decomposition ${\Sigma}={\mathcal H}\cup {\mathcal H}^*$ that if a non-trivial element
$h\in H_1({\Sigma})$ is not in the image of $\alpha$ then it intersects non-trivially with $H_1(\partial {\mathcal H})$, so
in this case $h\notin V(H)^{\perp}$. this implies that $V(H)^{\perp}\subset {\rm image}({\alpha})$, so $V(H)^{\perp}\subset
(V(H)\cap V(H)^{\perp})\cup V(H^*)=V(H^*)$. Therefore $V(H)^{\perp}= V(H^*)$, and this completes the proof of lemma
\ref{duality lemma 2}.

\qed

\bigskip

\section{The generalized Kauffman bracket and Jones polynomial of links on surfaces} \label{Jones section}

Various relations between the Tutte polynomial and link polynomials are well known, for example see \cite{Th, Ja}.
More recently \cite{ChP} such relations have been established for link polynomials and the Bollob\'{a}s-Riordan polynomial of
associated graphs on surfaces. (See also \cite{DFKLS, Mo} for other related results.)
In this section we consider a $2$-variable generalization of the Jones polynomial of links
in (surfaces)$\times I$, and more generally a $4-$variable Kauffman bracket of link diagrams on a surface, and we establish an analogue of
Thistlethwaite's theorem \cite{Th} relating these polynomials for alternating links
in ${\Sigma}\times I$ and the polynomial
$P_G$ of the associated Tait graph on the surface $\Sigma$. Using the interpretation of virtual links
as ``irreducible'' embeddings of classical links into surfaces \cite{Ku},
these results apply to virtual links.

Let $L$ be a link embedded in ${\Sigma}\times I$, where ${\Sigma}$ is a closed orientable surface.
Consider a projection $D$ of $L$ onto the surface. By general position $D$ is a diagram with a finite
number of crossings. Each crossing may be resolved as shown in figure \ref{resolutions}.
Given a diagram $D$ with $n$ crossings, consider the set ${\mathcal S}$ of its $2^n$ resolutions.
Each resolution $S\in {\mathcal S}$ is a disjoint collection of closed curves embedded in ${\Sigma}$. Denote by
${\alpha}(S)$ the number of resolutions of type (1) that were used to create it, by ${\beta}(S)$ the number of
resolutions of type (2), and let $c(S)$ be the number of
components of $S$. Consider the inclusion map $i\co S\subset {\Sigma}$,  and denote
$$k(S)\; =\; {\rm rank}\, (\, {\rm ker}\, \{i_*\co H_1(S)\longrightarrow H_1({\Sigma})\}).$$

\begin{figure}[ht]
\medskip
\includegraphics[height=2.2cm]{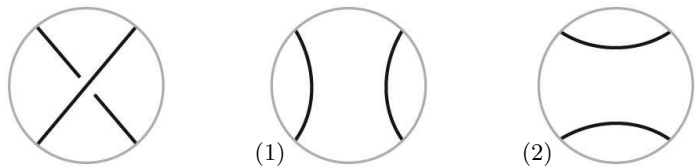}
{\scriptsize
    \put(-169,3){$(1)$}
    \put(-68,3){$(2)$}}
    \medskip
\caption{Resolutions of a crossing.}
\label{resolutions}
\end{figure}

Generalizing the classical definition of the Kauffman bracket (cf \cite{K1, B}), consider

\begin{equation} \label{multi Kauffman} \widetilde K_L(A,B,d) \; =\; \sum_{S\in {\mathcal S}} [i_*(H_1(S))]\, A^{{\alpha}(S)} \, B^{{\beta}(S)}\, d^{k(S)}
\end{equation}

Here $[i_*(H_1(S))]$
denotes a formal variable associated to the subgroup $i_*(H_1(S))$ of $H_1({\Sigma})$.
This definition is closely related to (more precisely, it may be viewed as a specialization of)
the {\em surface bracket polynomial},
defined in the context of virtual links in \cite{DK, M}.
Two diagrams in ${\Sigma}$, representing isotopic embeddings of a link $L$ in ${\Sigma}\times I$, are related by the usual
Reidemeister moves, and the usual specialization

\begin{equation} \label{Jones definition0}  \widetilde J_L(t)\; =\; (-1)^{w(L)}\, t^{3w(L)/4}\, \widetilde K_D (t^{-1/4}, t^{1/4},-t^{1/2} - t^{-1/2}),
\end{equation}

where $w(L)$ is the writhe of $L$, is an invariant of  an embedded oriented link $L\subset{\Sigma}\times I$. The polynomial $\widetilde J_L(t)$
with coefficients corresponding to subgroups of $H_1({\Sigma})$ may be used to distinguish non-isotopic
links in ${\Sigma}\times I$ (also see remark 6 following theorem \ref{Jones Tutte lemma} below.)
However if one is interested in studying links up to to the action of the diffeomorphisms of ${\Sigma}$, or in
studying virtual links, a relevant invariant is the following finite-variable specialization.
Denoting the rank of $i_*(H_1(S))$ by $r(S)$, define

\begin{equation} \label{Kauffman} K_D(A,B,d,Z) \; =\; \sum_{S\in{\mathcal S}} \, A^{{\alpha}(S)} \, B^{{\beta}(S)}\, d^{k(S)}\, Z^{r(S)},
\end{equation}
and the corresponding version of the Jones polynomial:
\begin{equation} \label{Jones definition}  J_L(t,Z)\; =\; (-1)^{w(L)}\, t^{3w(L)/4}\, K_D (t^{-1/4}, t^{1/4},-t^{1/2} - t^{-1/2},Z).
\end{equation}
Note that all of the polynomials considered here may be defined for {\em virtual} links \cite{K}, using their
``irreducible'' embeddings into surfaces \cite{Ku}.
Since $k(S)+r(S)$ equals the number $c(S)$ of components of $S$, it follows that for a virtual link
$L$, the invariant $K_D$ defined above specializes to the usual Kauffman bracket by setting $Z=d$:

$$ [L](A,B,d)\; =\; d^{-1}\, K_L(A,B,d,d).$$

\begin{figure}[ht]
\medskip
\includegraphics[height=2.2cm]{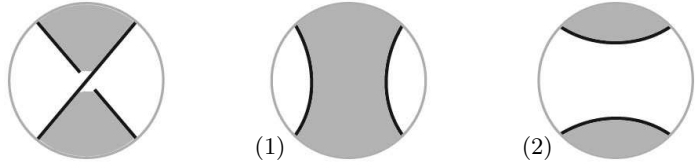}
{\scriptsize
    \put(-169,3){$(1)$}
    \put(-68,3){$(2)$}}
    \medskip
\caption{Checkerboard coloring near a crossing of an alternating diagram.}
\label{crossing}
\end{figure}

We now turn to the analogue for links on surfaces of Thistlethwaite's theorem \cite{Th} relating the Jones polynomial $J_L(t)$ of an alternating
link $L$ in $S^3$
to the specialization $T_G(-t, -t^{-1})$ of the Tutte polynomial of an associated Tait graph. Suppose $L$ is a link in ${\Sigma}\times I$ which has an
alternating diagram $D$ on $\Sigma$. Then this diagram may be checkerboard-colored, as shown near each crossing on the
left in figure \ref{crossing}. The associated Tait graph is the graph $G_D\subset {\Sigma}$ whose vertices correspond to
the shaded regions of the diagram, and two vertices are connected by an edge whenever the corresponding shaded regions
meet at a crossing (an example of an alternating link on the torus and the corresponding Tait graph are shown in figure \ref{alternating knot} -
compare with the example in \cite{ChP}.) The Tait graph is a well-defined graph $G\subset {\Sigma}$ if each component
in the complement of a link diagram $D$ is a disk; this condition holds for virtual links due to the irreducibility of
their embedding into ${\Sigma}\times I$.

\begin{theorem} \label{Jones Tutte lemma} \sl
The generalized Kauffman bracket (\ref{Kauffman}) of an alternating link diagram $D$
on a surface $\Sigma$ may be obtained from the polynomial $P_G$, defined by (\ref{poly def equation}),
of the associated Tait graph $G$ as follows:
\begin{equation} \label{Jones Tutte equation1}
K_D(A,B,d,Z)\; =\;
A^{g+v(G)-c(G)}\, B^{-g+n(G)}\, d^{c(G)}\, Z^{g}\,
P_{G}\left(\frac{Bd}{A}, \frac{Ad}{B},\frac{A}{BZ},\frac{B}{AZ}\right).
\end{equation}
In particular, substituting $A=t^{-1/4}, B=t^{1/4}, d=-t^{1/2}-t^{-1/2}$ as in (\ref{Jones definition}) yields an
expression for the $2$-variable Jones polynomials $J_L(t,Z)$ in terms of the polynomial $P_G$ of the associated Tait
graph.
\end{theorem}

\begin{figure}[t]
\bigskip
\includegraphics[height=3.5cm]{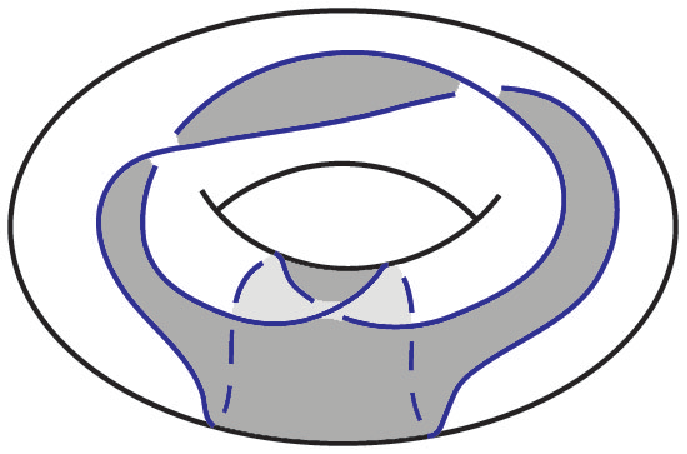} \hspace{1cm} \includegraphics[height=3.5cm]{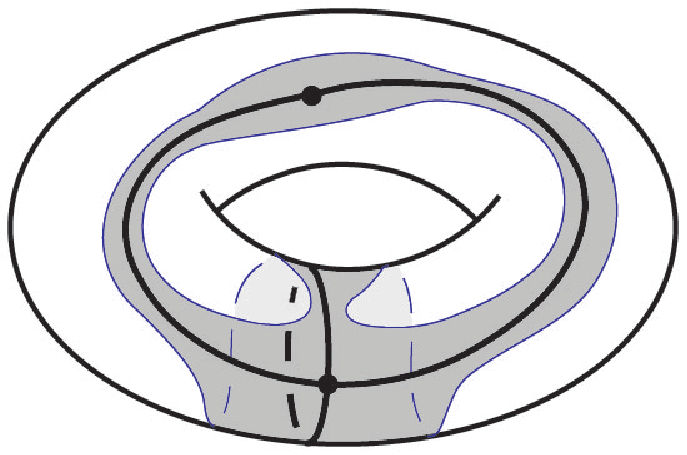}
\medskip
\caption{An alternating link diagram (left) and its Tait graph (right) on the torus.}
\label{alternating knot}
\end{figure}

{\bf Remarks.}\nl
{\bf 1.} Recall that $\widetilde K_L(A,B,d)$, defined in (\ref{multi Kauffman}), is a polynomial in $A,B,d$ with coefficients
corresponding to subgroups $V<H_1({\Sigma})$.
The equation (\ref{Jones Tutte equation1}) follows from a more general relation, which can be deduced from the proof of theorem \ref{Jones Tutte lemma}, between the polynomials
$\widetilde K_L$ and $\widetilde P_{G}(X,Y)$ (defined in (\ref{multi poly def equation})).
In particular, each coefficient
$[V]$, $V<H_1({\Sigma})$, of $\widetilde P_{G}$ is replaced with $V\cap V^{\perp}$ to get the corresponding coefficient
of $\widetilde K_L$.

{\bf 2.} It would be interesting to establish a relation, analogous to (\ref{Jones Tutte equation1})  between
these generalized versions of the Kauffman bracket, the Jones polynomial, and the Bollob\'{a}s-Riordan polynomial. In principle,
such a relationship follows from theorem \ref{Jones Tutte lemma} (see the discussion following theorem \ref{universality}), but an explicit formula does not seem to be
as straightforward as (\ref{Jones Tutte equation1}).

{\bf 3.} Theorem \ref{Jones Tutte lemma} asserts that the generalized Kauffman bracket $K_D(A,B,d,Z)$ may be obtained as
a specialization of the polynomial $P_G(D)$. It would be interesting to find out whether $K_D$ and $P_{G}$
(or $K_D$ and the Bollob\'{a}s-Riordan polynomial $BR_G$) in fact determine each other.

{\bf 4.} Suppose $D$ is an alternating link diagram (associated to a link $L\subset {\Sigma}\times I$) on an orientable surface ${\Sigma}$.
Switching each crossing,
one gets an alternating link diagram $D'$ whose checkerboard coloring is precisely that of $D$ with the colors switched on
each face. Therefore the associated graphs $G$, $G'$ are duals of each other. (To make this statement precise, it is
convenient to consider virtual links, so the embedding $L\subset {\Sigma}\times I$ is ``irreducible'' \cite{Ku}, and then the Tait graphs $G, G'$
are cellulations.) In this context theorem \ref{Jones Tutte lemma} gives a different perspective on the duality relation
(\ref{duality lemma equation}) for $P_G$.

{\bf 5.} Adapting the proof in \cite{ChP}, one may establish a generalization of theorem \ref{Jones Tutte lemma} from alternating diagrams
to checkerboard-colored diagrams, using a {\em signed} version of the polynomial $P_G$. The proof uses
the observation \cite{Ka} that any such link diagram on a surface can be made alternating by switching some of the
crossings, and then one labels by $-1$ each edge of the Tait graph
where a switch has been made.

{\bf 6.} One may generalize the correspondence between the Jones polynomial and $P_G$ to {\em skein modules.} (This is a further generalization from the polynomials $\widetilde K_L$ and $\widetilde P_G$ whose coefficients are subgroups of $H_1({\Sigma})$.)   Specifically,
one may consider the isotopy classes of
graphs on ${\Sigma}$ modulo the contraction-deletion relation, see section \ref{Tutte skein},
and the skein module of links modulo the Kauffman skein relation in figure \ref{resolutions}, cf. \cite{P, Tu}. The author would
like to thank J\'{o}zef Przytycki for pointing out this perspective on the problem.

\medskip

{\em Proof of theorem \ref{Jones Tutte lemma}.}
The terms in the expansions (\ref{poly def equation}), (\ref{Kauffman}) are in $1-1$ correspondence.
Specifically, for each spanning subgraph $H\subset G(D)$ parametrizing the sum (\ref{poly def equation}), consider the corresponding resolution $S(H)$: each crossing of the diagram $D$ is resolved as in figure \ref{crossing}, where the resolution (1) is used if the corresponding edge is included
in $H$, and the resolution (2) is used otherwise. Observing the effect of the resolutions on the shaded regions in figure \ref{crossing}, note that the collection of embedded curves $S(H)\subset {\Sigma}$ is the boundary of a regular neighborhood of $H$ in $\Sigma$. Moreover, the number ${\alpha}(S)$
of resolutions of type (1) is precisely the number $e(H)$ of edges of $H$, and ${\beta}(S)$ equals $e(G)-e(H)$.

Observe $${\alpha}(S(H))=e(H)=v(H)-c(H)+n(H),$$ $${\beta}(S(H))=e(G)-e(H)=n(G)-n(H)+c(H)-c(G).$$
Also note that since $S$ is the boundary of a regular neighborhood of $H$, $r(S)=l(H)$, and
$k(S)=c(H)+k(H)$. Therefore the summands $A^{{\alpha}(S)} \, B^{{\beta}(S)}\, d^{k(S)}\, Z^{r(S)}$ in (\ref{Kauffman}) may be
rewritten as $$ A^{v(H)-c(H)+n(H)} \,
B^{n(G)-n(H)+c(H)-c(G)}\, d^{c(H)+k(H)}\, Z^{l(H)}.$$
Substituting the required variables, the summands in the expansion (\ref{poly def equation}) of $P_G$ are of the form
$$\left(\frac{Bd}{A}\right)^{c(H)-c(G)} \left(\frac{Ad}{B}\right)^{k(H)} \left(\frac{A}{BZ}\right)^{s(H)/2}
\left(\frac{B}{AZ}\right)^{s^{\perp}(H)/2}.$$
The proof is completed by using the relations (\ref{identities}) to identify the exponents of $A,B,d,Z$ on the two sides of
(\ref{Jones Tutte equation1}).
\qed

\bigskip

\section{A multivariate graph polynomial} \label{multivariable section}

We conclude the paper by pointing out a multivariate version of the polynomial $P_G$, and observing
the corresponding duality relation. (Note that a multivariate version of the Bollob\'{a}s-Riordan polynomial has been considered in \cite{Mo}.
A duality result for a certain specialization of the signed multivariate Bollob\'{a}s-Riordan polynomial has been established in \cite{VT}.)
Let $G$ be a graph on a surface $\Sigma$, and let
$${\mathbf v}\, =\, \{ v_e\}_{e\in E(G)}$$
be a collection of commuting indeterminates associated to the edges of $G$. Following the notation used in (\ref{poly def equation}),
consider

\begin{equation} \label{multivariate surface Tutte}
\overline P_G(q,{\mathbf v}, A,B)\; =\; \sum_{H\subset G} q^{c(H)}\; A^{s(H)/2}\;
B^{s^{\perp}(H)/2} \! \prod_{e\in E(H)} v_e
\end{equation}

Clearly, the ``usual'' multivariate Tutte polynomial $Z_G$ \cite{S} is a specialization of $\overline P_G$:

$$ Z_G(q,{\mathbf v})\; = \; \overline P_{G,{\Sigma}}(q,{\mathbf v},1,1),$$

The relation to the polynomial $P_G(X,Y,A,B)$ defined in (\ref{poly def equation}) is given by

$$ P_G(X,Y,A,B)\; = \; X^{-c(G)} \, Y^{-g-v(G)}\, \overline P_G(XY,\,Y,\,A/Y,\,BY),$$

where as usual $c(G)$ denotes the number of connected components of the graph $G$, $v(G)$ is the number
of vertices of $G$, and $g$ is the genus of the surface $\Sigma$.
That is, to get the polynomial $P_G$, one sets in the multivariate version $\overline P_G$ all edge
weights
$v_e$ equal to $Y$, and $q=XY$.
The analogue of the duality (\ref{duality lemma equation}) for the multivariate polynomial $\overline P_G$ is as follows.

\begin{lemma} \label{mutlivariate duality lemma} Let $G$ be a cellulation of a surface $\Sigma$ (or equivalently a ribbon graph), and let $G^*$ denote its
dual. Then
\begin{equation} \label{multivariate duality}
\overline P_{G^*}(q,{\mathbf v}, A,B)\; =\; q^{-g+c(G^*)-v(G)}(\!\prod_{e\in E(G)} v_e)\; \overline P_G(q,q/{\mathbf v}, B/q,Aq).
\end{equation}
\end{lemma}

As the notation indicates, the edge weights of $G^*$ in the formula on the right-hand side
are given by $\{q/v_e\}_{e\in E}$. Using the relation $c(H)=v(H)-e(H)+n(H)$, note that the expansion of the polynomial
$\overline P_G$ may be rewritten as

\begin{equation}
\overline P_G(q,{\mathbf v}, A,B)\; =\; q^{v(G)} \sum_{H\subset G} q^{n(H)}\; A^{s(H)/2}\;
B^{s^{\perp}(H)/2} \! \prod_{e\in E(H)} \frac{v_e}{q}
\end{equation}
The proof of lemma \ref{mutlivariate duality lemma} consists of identifying the terms in the
expansions of the two sides, following the lines of the proof of theorem \ref{duality lemma}.

Note that the usual duality relation for {\em planar} graphs (cf. \cite{S}):
$$Z_{G^*}(q,{\mathbf v})\, =\, q^{1-v(G)}\left(\prod_{e\in E(G)}v_e \right) Z_G(q,q/{\mathbf v})$$
may be obtained as a specialization of the relation (\ref{multivariate duality}).

\end{document}